\documentclass[2p]{article}
\usepackage{amssymb,amsmath,amsthm}
\usepackage{indentfirst}
\usepackage{exscale}
\usepackage{relsize}
\usepackage[numbers,sort&compress]{natbib}

\usepackage{geometry}
\usepackage{color}
\newcommand{\R}{\mathbb{R}}
\theoremstyle{definition}

\allowdisplaybreaks
\theoremstyle{remark}

\numberwithin{equation}{section}
\UseRawInputEncoding

\begin{document}
	\title{\Large\bf{ Least energy solutions of two asymptotically cubic Kirchhoff equations on locally finite graphs}}
	\date{}
	\author {Zhangyi Yu$^{1}$, \ Xingyong Zhang$^{1,2}$\footnote{Corresponding author, E-mail address: zhangxingyong1@163.com} , \ Xin Ou$^{1}$\\
		{\footnotesize $^1$Faculty of Science, Kunming University of Science and Technology,}\\
		{\footnotesize Kunming, Yunnan, 650500, P.R. China.}\\
		{\footnotesize $^{2}$Research Center for Mathematics and Interdisciplinary Sciences, Kunming University of Science and Technology,}\\
		{\footnotesize Kunming, Yunnan, 650500, P.R. China.}\\}
	\date{}
	\maketitle

	\begin{center}
		\begin{minipage}{15cm}
			\par
			\small  {\bf Abstract:} We study the existence of least energy solutions for two Kirchhoff equations with the asymptotically cubic nonlinearity $f(u)=\lambda u+\eta|u|^2u$ on a locally weighted and connected finite graph $G=(V,E)$. Such nonlinearity satisfies neither $\frac{F(u)}{u^4}\to +\infty$ as $|u|\to\infty$, where $F(u)=\int_0^uf(s)ds$, nor $\frac{f(u)}{u}\to 0$ as $u\to 0$. By utilizing the  constrained variational method, we prove that there exist $\lambda_1\ge 0$ and  $\eta_0\ge 0$ ($\lambda_1^*\ge 0$ and $\eta_0^*\ge 0$) such that these two equations have at least a least energy  solution if $|\lambda|<a\lambda_1$ ($|\lambda|<a\lambda_1^*$) and $\eta>\eta_0$ ($\eta>\eta_0^*$).

		\end{minipage}
	\end{center}
	\allowdisplaybreaks
	\vskip2mm
	{\section{Introduction }}
	\setcounter{equation}{0}
	In this paper, we are concerned with the least energy solutions for the following two Kirchhoff equations  on a locally weighted and connected finite graph $G=(V,E)$:
	\begin{eqnarray}
		\label{c1}
		\begin{cases}
			-\left(a+b\int_{\Omega}|\nabla u|^2d\mu\right)\Delta u=\lambda u+\eta|u|^2u,\;&x\in \Omega,\\
			u=0,\;&x\in\partial \Omega,
		\end{cases}
	\end{eqnarray}
	and
	\begin{eqnarray}\label{t2}
		-\left(a+b\int_{V}(|\nabla u|^2+h(x) u^2)d\mu\right)(\Delta u+h(x)u)=\lambda u+\eta|u|^2u,\;x\in V,
	\end{eqnarray}
	where $V$ denotes the vertexes sets, $E$ denotes the edges sets, $\Omega\cup\partial\Omega\subset V$ is a bounded domain, $a>0$, $b\ge0$, $\lambda,\eta\in\R$, and $h:V\to \R$, $\nabla$ and $\Delta$ are the gradient operator and the Laplacian operator defined on the locally finite graph, respectively.
Next, we recall in details these basic concepts on locally finite graphs which can be seen in \cite{Yamabe 2016, Grigor2017, Grigor-book, Yang}. $G=(V, E)$ is called a locally connected finite graph if for any two vertexes in $V$ can be connected by finite edges in $E$ and for any $x\in V$ there are only finite edges $xy\in E$. The weight on any edge $xy\in E$ is represented by $\omega_{xy}$, which is required to satisfy $\omega_{xy}=\omega_{yx}$ and $\omega_{xy}>0$. We use $y\thicksim x$ to represent those vertices $y$ that connect with $x$. The distance  between two vertices $x$ and $y$, denoted by $d(x,y)$, is defined by the minimal number of edges which link $x$ to $y$. Assume $\Omega\subset V$. If $d(x,y)$ is uniformly bounded for any $x,y\in\Omega$, then $\Omega$ is known as a bounded domain in $V$. Let
	$$
	\partial\Omega=\{y\in V,\;y\notin\Omega|\;\exists\;x\in\Omega\;\text{such\;that}\;xy\in E\},
	$$
	which is known as the boundary of $\Omega$. Set $\Omega^{\circ}=\Omega\backslash\partial\Omega$ which is known as the interior of $\Omega$. Obviously, $\Omega^{\circ}=\Omega$.
	\par
	Assume that $\mu:V\rightarrow \R^+$ is a finite measure  and there exists a $\mu_0>0$ such that $\mu(x)\ge \mu_0$.
	The Laplacian operator $\Delta$ of $\psi$ is defined as
	\begin{eqnarray*}
		\label{eq3}
		\Delta \psi(x)=\frac{1}{\mu(x)}\sum\limits_{y\thicksim x}w_{xy}(\psi(y)-\psi(x)),\ \ \mbox{for all }x\in V.
	\end{eqnarray*}
	The corresponding gradient form has the expression given below
	\begin{eqnarray*}
		\label{eq4}
		\Gamma(\psi_1,\psi_2)(x)=\frac{1}{2\mu(x)}\sum\limits_{y\thicksim x}w_{xy}(\psi_1(y)-\psi_1(x))(\psi_2(y)-\psi_2(x)):=\nabla\psi_1\cdot\nabla\psi_2,\ \ \mbox{for all }x\in V.
	\end{eqnarray*}
	The length of the gradient is denoted by
	\begin{eqnarray*}
		\label{eq5}
		|\nabla \psi|(x)=\sqrt{\Gamma(\psi,\psi)(x)}=\left(\frac{1}{2\mu(x)}\sum\limits_{y\thicksim x}w_{xy}(\psi(y)-\psi(x))^2\right)^{\frac{1}{2}}, \ \ \mbox{for all }x\in V.
	\end{eqnarray*}
	For every fixed function $\psi:\Omega\rightarrow\mathbb{R}$ or $\psi:V\rightarrow\mathbb{R}$, we set
	\begin{eqnarray*}
		\label{eq7}
		\int_\Omega \psi(x) d\mu=\sum\limits_{x\in \Omega}\mu(x)\psi(x),\ \ \mbox{and  }  \int_V\psi(x) d\mu=\sum\limits_{x\in V}\mu(x)\psi(x)£¬
	\end{eqnarray*}
 respectively.
	\par
  To investigate the equation (\ref{c1}),  we denote  $C_0(\Omega)=\{\psi:\Omega\to \R|\psi=0 \mbox{ on } \partial \Omega\}$ and let
	$$
	C_c(\Omega)=\{\phi:V\rightarrow\R|supp\  \phi\subset\Omega\;\text{and}\;\phi(x)=0 \mbox{ for } \forall x\in V\backslash\Omega\}.
	$$
	For any function $\phi\in C_c(\Omega)$, there holds
	\begin{eqnarray*}
		\label{eq15}
		\int_{\Omega}\phi \Delta \psi d\mu=-\int_{\Omega\cup\partial\Omega}\Gamma(\psi,\phi) d\mu.
	\end{eqnarray*}
   For any given $1\le \gamma<+\infty$, let $L^\gamma(\Omega)$ be the completion of $C_c(\Omega)$ under the norm
	$$
	\|u\|_{\gamma}=\left(\int_{\Omega}|u(x)|^\gamma d\mu\right)^{\frac{1}{\gamma}}.
	$$
	Let $W_0^{1,2}(\Omega)$ be the completion of $C_c(\Omega)$ under the norm
	$$
	\|u\|=\left(\int_{\Omega\cup\partial\Omega}|\nabla u(x)|^2d\mu\right)^{\frac{1}{2}}.
	$$
	For any $u\in W_0^{1,2}(\Omega)$, we also define the following norm:
	$$
	\|u\|_{\infty,\Omega}=\max\limits_{x\in\Omega}|u(x)|.
	$$
	It is easy to see that $W_0^{1,2}(\Omega)$ is of finite dimension with  $\mbox{dim }W_0^{1,2}(\Omega)=\#\Omega$ where $\#\Omega$ denotes the number of elements in $\Omega$. We shall work on $W_0^{1,2}(\Omega)$ to investigate the equation (\ref{c1}).
	\par
	To investigate the equation (\ref{t2}), we denote  $C(V) =\{u|u:V\to \R\}$. For $u\in C(V)$, its support set is defined as $\text{supp}(u)=\{x\in V:u(x)\not=0\}$. Let $C_c(V)$ be the set of all functions with finite support, and for any function $\phi\in C_c(V)$, there holds
	\begin{eqnarray*}
		\label{eq15}
		\int_{V}\phi \Delta \psi d\mu=-\int_{V}\Gamma(\psi,\phi) d\mu.
	\end{eqnarray*}
 For any given $1\le \gamma<+\infty$, let $L^\gamma(\Omega)$ be the completion of $C_c(\Omega)$ under the norm
	$$
	\|u\|_{V,\gamma}=\left(\int_{V}|u(x)|^\gamma d\mu\right)^{\frac{1}{\gamma}}.
	$$
\par
 We make the following assumptions for $h$:\\
  {\it {\bf$(H_1)$} there exists a constant $h_0>0$ such that $h(x)\ge h_0>0$ for all $x\in V$;\\
  {\bf $(H_2)$}  $h(x)\to +\infty$ as $d(x,x_0)\to\infty$ for some fixed $x_0$.}

  \vskip2mm
  \par
 Let $W^{1,2}(V)$ be the completion of $C_c(V)$ under the norm
	$$
	\|u\|_{W^{1,2}(V)}=\left(\int_V\left(|\nabla u(x)|^2+|u(x)|^2\right)d\mu\right)^\frac{1}{2}.
	$$
	Define the space
	$$
	W_h^{1,2}(V)=\left\{u\in W^{1,2}(V) \Big|\int_V h(x)|u(x)|^2 d\mu<\infty\right\}
	$$
	endowed with the norm
	$$
	\|u\|_{W_h^{1,2}(V)}=\left(\int_V\left(|\nabla u(x)|^2+h(x)|u(x)|^2\right)d\mu\right)^\frac{1}{2}.
	$$
 $W_h^{1,2}(V)$ is a reflexive Banach space.	We shall work on $W_h^{1,2}(\Omega)$ to investigate the equation (\ref{t2}).
 \par
	In recent years, the investigations on solutions for the equations on locally  finite graphs have attracted lots of attention because of their applications to image processing, data analysis and the unique mathematical properties of graphs, for example, see \cite{Elmoataz2012,Elmoataz2015,Elmoataz2017, Shao2023,Han2020,Han2021,Grigor2017,Pinamonti2022,Zhang2018,Zhang2019,Zhang 2022,Pang 2024,Lin2017,Liu2024,Qiu2023,Yang 2023,Yang,Xu2023,Yang2024}.
	Especially, in \cite{Pan2023}, Pan-Ji studied the least energy sign-changing solutions to the following nonlinear Kirchhoff equation
	\begin{eqnarray*}
		\label{c2}
		-\left(a+b\int_V|\nabla u|^2d\mu\right)\Delta u+c(x)u=f(u)
	\end{eqnarray*}
	on a locally finite graph $G=(V,E)$, where $a,b$ are positive constants and $c:V\to \R$. They used the constrained variational method to prove the existence of a least energy sign-changing solution when $f$ is assumed to satisfy the superlinear conditions.
	They also showed the sign-changing least energy is strictly larger than twice that of the least energy.
 In \cite{Ou2024}, Ou-Zhang also obtained the similar result  for the following Kirchhoff-type equations with general power law, logarithmic nonlinearity and Dirichlet boundary value on a locally finite graph:
	\begin{eqnarray*}
\label{eq1}
 \begin{cases}
  -\left(a+b\int_{\overline{\Omega}}|\nabla u|^2d\mu\right)\Delta u+\lambda g(x)u^{\frac{2k}{m}-1}=Q(x)|u|^{p-2}u\ln |u|^r,& \text {in} \; \Omega^\circ,\\
  u=0,&\text {on} \; \partial\Omega,\\
   \end{cases}
\end{eqnarray*}
where $\overline{\Omega} \subset V $ is a bounded domain, $\overline{\Omega}:=\Omega^\circ\cup\partial \Omega=\Omega\cup\partial \Omega$, $a>0$, $b\geq0$, $p>4$, $\lambda\ge0$, $r\geq 1$, $m,k\in \mathbb{N}$ satisfying $1<\frac{2k}{m}\leq p$ and $Q,g\in\Omega^\circ\to \R^+$.
In both \cite{Ou2024} and \cite{Pan2023}, the nonlinear terms $f$ are required to satisfy the following superlinear condition:\\
{\it $(C_1)$ $\frac{F(u)}{u^4}\to +\infty$, as $|u|\to\infty$, where $F(u)=\int_0^uf(s)ds$.}\\
{\it $(C_2)$ $\frac{f(u)}{u}\to 0$, as $u\to 0$.}\\
\vskip2mm
\par
Recently, in \cite{Li2024}, Li-Wang considered the following equation
\begin{eqnarray}
		\label{pp1}
		-\Delta u+(\lambda a(x)+1)u=f(u), \ \ \forall x\in V,
	\end{eqnarray}
where $a:V\to [0,+\infty)$, $\lambda>1$ and $f$ is required to satisfy an asymptotically linear condition: $\frac{f(u)}{u}\to l$ for some constant $l$ with a non-negative lower bound. However, they also needed $(C_2)$ to hold. They proved the existence of a least energy solution $u_{\lambda}$  for (\ref{pp1}) and the convergence of least energy solution family $\{u_{\lambda}\}$ as $\lambda\to +\infty$.

\par
A natural question is whether the least energy solution still exists if both the conditions $(C_1)$ and $(C_2)$ are unavailable. In this paper, we establish the existence of the least energy solutions for (\ref{c1}) and (\ref{t2}).  It is easy to see that the nonlinear term $\lambda u + \eta|u|^2u$ in (\ref{c1}) and (\ref{t2})  is neither satisfying $(C_1)$ nor $(C_2)$. Our work is also motivated by \cite{Khoutir} and \cite{Li2022} where the authors considered the existence of sign-changing least energy solutions for Shr\"odinger-possion and a Kirchhoff equation with  the nonlinear term $\lambda u + \eta|u|^2u$ in Euclidean setting. We develop the idea in \cite{Khoutir} and \cite{Li2022} to the locally finite graph setting, but we do not consider the existence of sign-changing least energy solutions. That is because it involves more complex factors in locally finite graph setting than in the Euclidean setting when we decompose $\|u\|$ into $\|u^+\|+\|u^-\|$ (for example, see \cite{Ou2024}) so that it is not easy to obtain a result similar to Lemma 2.2 in \cite{Li2022} or Lemma II.2 in \cite{Khoutir}. However, the result holds when we consider the least energy solutions (see Lemma 2.6 and Lemma 3.5 below). Next, we state our main results.
\vskip2mm
	\par
	Let
	\begin{eqnarray}\label{oo1}
	\lambda_1=\inf_{u\in W^{1,2}_0(\Omega)\setminus\{0\}}\frac{\|u\|^2}{\int_{\Omega}|u|^2d\mu},\ \  \eta_0=\inf_{u\in W^{1,2}_0(\Omega)\setminus\{0\}}\frac{b\|u\|^4}{\int_{\Omega}|u|^4d\mu}
	\end{eqnarray}
and
 \begin{eqnarray}\label{oo2}
 \lambda_1^*=\inf_{u\in W^{1,2}_h(V)\setminus\{0\}}\frac{\|u\|^2_{W^{1,2}_h(V)}}{\int_{V}|u|^2d\mu},\ \  \eta_0^*=\inf_{u\in W^{1,2}_h(V)\setminus\{0\}}\frac{b\|u\|^4_{W^{1,2}_h(V)}}{\int_{V}|u|^4d\mu}.
 \end{eqnarray}

	\vskip2mm
	\noindent
	{\bf Remark 1.1} Obviously, $\lambda_1,\lambda_1^*,\eta_0,\eta_0^*\ge 0$, and if $b=0$, then $\eta_0=\eta_0^*=0$. Furthermore, a simple argument can show that $\lambda_1>0$. Indeed, assume that $\lambda_1=0$. Then there exists a sequence $\{u_n\}\subset W_0^{1,2}(\Omega)\setminus\{0\}$ such that $\frac{\|u_n\|^2}{\int_{\Omega}|u_n|^2d\mu}\to 0$ as $n\to \infty$. On the other hand, since $W_0^{1,2}(\Omega)$ is of finite dimension, the norm $\|u_n\|$ is equivalent to $\|u_n\|_2$. Hence, there exists a constant $D_0>0$ such that $\frac{\|u_n\|^2}{\int_{\Omega}|u_n|^2d\mu}\ge D_0$ for all $n\in\mathbb N$. A contradiction appears. Hence, $\lambda_1>0$. Similarly, by using the equivalence of $\|u_n\|$ and $\|u_n\|_4$, we can also prove that $\eta_0>0$ if $b>0$. Moreover, for all $\{u_n\}\subset W_h^{1,2}(V)$, note that $\|u_n\|_{W^{1,2}_h(V)}\ge h_0^{\frac{1}{2}}\|u_n\|_{V,2}$ and $\|u_n\|_{W^{1,2}_h(V)}\ge \mu_0^{\frac{1}{4}} h_0^{\frac{1}{2}}\|u_n\|_{V,4}$ (see Lemma 3.1 below). Similar to the above arguments, we can also prove that $\lambda_1^*>0$ and if $b\not=0$, then $\eta_0^*>0$.
	\vskip2mm
	\noindent
	{\bf Theorem 1.1.} {\it If $\eta\le\eta_0$, then $(\ref{c1})$ has no any nontrivial solution. If $|\lambda|<a\lambda_1$ and $\eta>\eta_0$, then $(\ref{c1})$ possesses a least energy  solution $u_b$ and $\|u_b\|>\left(\frac{a-\frac{\lambda}{\lambda_1}}{\eta D_{\Omega}^4}\right)^{\frac{1}{2}}$ and
		$$
		I(u_b)>\frac{1}{4}\left(a-\frac{|\lambda|}{\lambda_1}\right)\cdot\left(\frac{a-\frac{\lambda}{\lambda_1}}{\eta D_{\Omega}^4}\right),
		$$
		where $D_{\Omega}$ is the embedding constant from $W_0^{1,2}(\Omega)$ to $L^4(\Omega)$ and $I(\cdot)$ is the corresponding energy functional of equation (\ref{c1}) and is given by}
\begin{eqnarray}\label{cb2}
		I(u)=\frac{a}{2}\|u\|^2-\frac{\lambda}{2}\int_{\Omega}|u|^2d\mu+\frac{b}{4}\|u\|^4-\frac{\eta}{4}\int_{\Omega}|u|^4d\mu,\;\;\;\forall u\in W^{1,2}_0(\Omega).
	\end{eqnarray}

	\vskip2mm
	\noindent
	{\bf Theorem 1.2.} {\it Assume that $(H_1)$ and $(H_2)$ hold. If $\eta\le\eta_0^*$, then $(\ref{t2})$ has no any nontrivial solution. If $|\lambda|<a\lambda_1^*$ and $\eta>\eta_0^*$, then $(\ref{t2})$ possesses a least energy  solution $u_b^*$ and $\|u_b^{*}\|>\left(\frac{a-\frac{\lambda}{\lambda_1^*}}{\eta}\mu_0h_0^2\right)^{\frac{1}{2}}$ and
		$$
		J(u_b^*)>\frac{1}{4}\left(a-\frac{|\lambda|}{\lambda_1^*}\right)\cdot\left(\frac{a-\frac{\lambda}{\lambda_1^*}}{\eta}\mu_0h_0^2\right),
		$$
		where $J(\cdot)$ is the corresponding energy functional of equation (\ref{t2}) and is given by}
\begin{eqnarray}\label{cm2}
		J(u)=\frac{a}{2}\|u\|_{ W^{1,2}_h(V)}^2-\frac{\lambda}{2}\int_{V}|u|^2d\mu+\frac{b}{4}\|u\|_{ W^{1,2}_h(V)}^4-\frac{\eta}{4}\int_{V}|u|^4d\mu,\;\;\forall u\in W^{1,2}_h(V).
	\end{eqnarray}

	\par
When $a=1$ and $b=0$,  (\ref{c1}) and (\ref{t2}) reduce to the following equations, respectively:
	\begin{eqnarray}
		\label{k1}
		\begin{cases}
			-\Delta u=\lambda u+\eta|u|^2u,\;&x\in \Omega,\\
			u=0,\;&x\in\partial \Omega.
		\end{cases}
	\end{eqnarray}
	and
	\begin{eqnarray}
		\label{kk1}
		-\Delta u+h(x)u=\lambda u+\eta|u|^2u,\;&x\in V.
	\end{eqnarray}
	Then Theorem 1.1 and Theorem 1.2 reduce to the following results correspondingly.
	
	\vskip2mm
	\noindent
	{\bf Corollary 1.1.} {\it If $\eta\le0$, then $(\ref{c1})$ has no any nontrivial solution. If $|\lambda|<\lambda_1$ and $\eta>0$, then $(\ref{c1})$ possesses a least energy  solution $u_0^*$ and $\|u_0^{*}\|>\left(\frac{1-\frac{\lambda}{\lambda_1}}{\eta  D_{\Omega}^4}\right)^{\frac{1}{2}}$ and}
		$$
		I(u_0^*)>\frac{1}{4}\left(1-\frac{|\lambda|}{\lambda_1}\right)\cdot\left(\frac{1-\frac{\lambda}{\lambda_1}}{\eta  D_{\Omega}^4}\right).
		$$
	
\vskip2mm
	\noindent
	{\bf Corollary 1.2.} {\it If $\eta\le0$, then $(\ref{t2})$ has no any nontrivial solution. If $|\lambda|<\lambda_1^*$ and $\eta>0$, then $(\ref{t2})$ possesses a least energy  solution $u_0^*$ and $\|u_0^{*}\|>\left(\frac{1-\frac{\lambda}{\lambda_1^*}}{\eta}\mu_0h_0^2\right)^{\frac{1}{2}}$ and}
		$$
		J(u_0^*)>\frac{1}{4}\left(1-\frac{|\lambda|}{\lambda_1^*}\right)\cdot\left(\frac{1-\frac{\lambda}{\lambda_1^*}}{\eta}\mu_0h_0^2\right).
		$$

	\vskip0mm
	\noindent
	{\bf Remark 1.2.} We remark that in \cite{Li2022} and \cite{Khoutir}, $\lambda$ is required to satisfy a weaker condition that $\lambda<a\lambda_1$. However, by checking the proofs in \cite{Li2022} and \cite{Khoutir}, it is not difficult to see that in \cite{Li2022} and \cite{Khoutir}, $\lambda$ actually has to satisfy $0\le\lambda<a\lambda_1$ (see the proofs of Lemma 2.3 and Lemma 2.4 in \cite{Li2022} or Lemma II.2 and Lemma II.3 in \cite{Khoutir}). In our proofs, the restriction can be relaxed to $|\lambda|<a\lambda_1$ by changing $\lambda$ to  $|\lambda|$ simply (see the proofs of Lemma 3.2 and Lemma 3.3 below).
	
  \vskip0mm
	\noindent
	{\bf Remark 1.3.}  Corollary 1.1 and Corollary 1.2 are different from those results in \cite{Li2024} since the nonlinearity $\lambda u+\eta|u|^2u$ fails to satisfy the condition $(C_2)$ and is super-linear at infinity rather than  asymptotically linear for equations (\ref{k1}) and (\ref{kk1}).

	\vskip2mm
	{\section{The least energy solution for $(\ref{c1})$}}
	\setcounter{equation}{0}
	
	In this section, we prove the existence of the least energy solution  for $(\ref{c1})$ and complete the proof of Theorem 1.1. We shall work on $W_0^{1,2}(\Omega)$ and need  the following Sobolev embedding theorem.
\vskip2mm
	\noindent
	{\bf Lemma 2.1.}  (\cite{Yamabe 2016}, Theorem 7 with $m=1$ and $l=2$) {\it Suppose that $G=(V,E)$ is a locally finite graph, $\Omega$ is a bounded domain of $V$ such that $\Omega^{\circ}\neq\emptyset$. Then $W_0^{1,2}(\Omega)$ is embedded in $L^{\gamma}(\Omega)$ for all $1\leq \gamma\leq+\infty$. Particularly, there exists a constant $D_{\Omega}>0$ depending only on  $\Omega$ such that
		\begin{eqnarray*}
			\left(\int_{\Omega}|u|^\gamma d\mu\right)^{\frac{1}{\gamma}}\leq D_\Omega\left(\int_{\Omega\cup\partial\Omega}|\nabla u|^2d\mu\right)^{\frac{1}{2}}, \ \  \|u\|_{\infty,\Omega}\le M\|u\|_{W_0^{1,2}(\Omega)}
		\end{eqnarray*}
		for all $1\le \gamma<+\infty$ and all $u\in W_0^{1,2}(\Omega)$, and $D_\Omega=\frac{C}{\mu_{\min,\Omega}}(1+\sum_{x\in\Omega}\mu(x))$ with $C$ satisfying $\|u\|_{L^2(\Omega)}\le C\|u\|_{W_0^{1,2}(\Omega)}$, $M=\frac{C}{\mu_{\min,\Omega}^\frac{1}{2}}$ and $\mu_{\min,\Omega}=\min_{x\in\Omega}\mu(x)$. Moreover, $W_0^{1,2}(\Omega)$ is pre-compact, that is, if $\{u_n\}$ is bounded in $W_0^{1,2}(\Omega)$, then up to a subsequence, there exists some $u\in W_0^{1,2}(\Omega)$ such that $u_n\rightarrow u$ in $W_0^{1,2}(\Omega)$.}
	\vskip2mm
	\par
	Let $I_{\lambda}(u):=a\|u\|^2-\lambda\int_{\Omega}|u|^2d\mu$ and $I_{\eta}(u):=b\|u\|^4-\eta\int_{\Omega}|u|^4d\mu$. Then the energy functional given by (\ref{cb2}) can be written as
\begin{eqnarray}\label{b2}		
I(u)=\frac{1}{2}I_{\lambda}(u)+\frac{1}{4}I_{\eta}(u),\;\;\;\forall u\in W^{1,2}_0(\Omega).
	\end{eqnarray}
Moreover, a standard argument implies that $I(u)$ belongs to $C^1(W_0^{1,2}(\Omega),\R)$, and
	\begin{eqnarray}\label{3.2}
		\langle I'(u),v\rangle=a\int_{\Omega\cup\partial\Omega}\nabla u\nabla vd\mu-\lambda\int_{\Omega}uvd\mu+b\|u\|^2\int_{\Omega\cup\partial\Omega}\nabla u\nabla vd\mu-\eta\int_{\Omega}|u|^2uvd\mu
	\end{eqnarray}
	for any $u,v\in W^{1,2}_0(\Omega)$ (for example, see \cite{Yang}). Thus, the critical points of $I$ correspond to the weak solutions of (\ref{c1}), which actually are the point-wise solutions of (\ref{c1}) (for example, see  Proposition 3.1 in \cite{Ou2024}).
	
	\par
	We denote the following constraint
	\begin{eqnarray}\label{oo4}
	\mathcal{K}=\{u\in W^{1,2}_0(\Omega),u\not=0:\langle I'(u),u\rangle=0\},
	\end{eqnarray}
	and
	$$
	c=\inf_{u\in \mathcal{K}}I(u).
	$$
	If $u\in W^{1,2}_0(\Omega)$ is a nontrivial solution of $(\ref{c1})$, then by (\ref{3.2}), it holds that
	\begin{eqnarray}
		\label{b3}
		0=\langle I'(u),u\rangle=I_{\lambda}(u)+I_{\eta}(u).
	\end{eqnarray}
	Thus, we deduce from $\lambda<a\lambda_1$, $(\ref{b3})$ and $(\ref{oo1})$ that
	\begin{eqnarray*}
		I_{\eta}(u)=b\|u\|^4-\eta\int_{\Omega}|u|^4d\mu= -I_{\lambda}(u)
		=-a\|u\|^2+\lambda\int_{\Omega}|u|^2d\mu<-a\|u\|^2+a\lambda_1\int_{\Omega}|u|^2d\mu\le0,
	\end{eqnarray*}
	that is
	\begin{eqnarray*}
		\eta>\frac{b\|u\|^4}{\int_{\Omega}|u|^4d\mu},\;\;\;\text{for any }u\in\mathcal{K}.
	\end{eqnarray*}
	\vskip2mm
	\noindent
	{\bf Lemma 2.2.} {\it Suppose $\eta\le\eta_0=\inf_{u\in W^{1,2}_0(\Omega)\setminus\{0\}}\frac{b\|u\|^4}{\int_{\Omega}|u|^4d\mu}$. Then $(\ref{c1})$ has no any nontrivial solution.}
	\vskip2mm
	\noindent
	{\bf Proof.} It follows from $\eta\le\eta_0$ that  for any $u\in W^{1,2}_0(\Omega)$ with $u\not =0$,
	\begin{eqnarray}\label{b4}
		I_\eta(u)=b\|u\|^4-\eta\int_{\Omega}|u|^4d\mu\ge 0.
	\end{eqnarray}
	Since $I_{\lambda}(u)>0$ by $\lambda<a\lambda_1$, then (\ref{b4}) implies that
	$$
	\langle I'(u),u\rangle=I_{\lambda}(u)+I_{\eta}(u)>0,
	$$
	which means that $\mathcal{K}=\emptyset$. The proof is completed. \qed
	\vskip2mm
	\par
	We denote
	\begin{eqnarray*}
		\mathcal{U}=\Bigg\{u\in W^{1,2}_0(\Omega),u\not=0:I_{\eta}(u)<0\Bigg\}
		=\left\{u\in W^{1,2}_0(\Omega),u\not=0:b\|u\|^4-\eta\int_{\Omega}|u|^4d\mu<0\right\}.
	\end{eqnarray*}

	\vskip2mm
	\noindent
	{\bf Lemma 2.3.} {\it If $\lambda<a\lambda_1$ and $\eta>\eta_0$, then $\mathcal{U}\not=\emptyset$ and $\mathcal{K}\subset\mathcal{U}$.}
	\vskip0mm
	\noindent
	{\bf Proof.}  Assume that $\eta>\eta_0$. By the definition of $\eta_0$, there exists a $v\in W^{1,2}_0(\Omega)$ such that $v\not=0$ and
	$$
	\eta>\frac{b\|v\|^4}{\int_{\Omega}|v|^4d\mu}\ge \eta_0.
	$$
	Thus, $b\|v\|^4-\eta\int_{\Omega}|v|^4d\mu<0$, which implies $v\in\mathcal{U}$. Hence, $\mathcal{U}\not=\emptyset$.
	Next, we prove $\mathcal{K}\subset\mathcal{U}$. For any $u\in\mathcal{K}$, we have
	\begin{eqnarray}
		\label{b5}
		b\|u\|^4-\eta\int_{\Omega}|u|^4d\mu=-I_{\lambda}(u).
	\end{eqnarray}
	Furthermore, by the fact that  $I_{\lambda}(u)>0$,  $(\ref{b5})$ shows that $u\in\mathcal{U}$. So $\mathcal{K}\subset\mathcal{U}$.
	\qed
	\vskip2mm
	\noindent
	{\bf Lemma 2.4.} {\it If $|\lambda|<a\lambda_1$ and $\eta>\eta_0$, then there exists $\kappa_{\lambda,\eta}>0$ such that $\|u\|>\kappa_{\lambda,\eta}$ for every $u\in\mathcal{K}$.
	}
	\vskip0mm
	\noindent
	{\bf Proof.}  For any $u\in\mathcal{K}$, using $(\ref{b3})$, $(\ref{oo1})$ and Lemma 2.1, we get
	\begin{eqnarray*}
		\left(a-\frac{|\lambda|}{\lambda_1}\right)\|u\|^2&\le& a\|u\|^2-|\lambda|\int_{\Omega}|u|^2d\mu\le I_{\lambda}(u)
		<  I_{\lambda}(u)+b\|u\|^4
		=
		\eta\int_{\Omega}|u|^4d\mu\le\eta D_{\Omega}^4\|u\|^4.
	\end{eqnarray*}
	 Therefore, $\|u\|>\left(\frac{a-\frac{|\lambda|}{\lambda_1}}{\eta D_{\Omega}^4}\right)^{\frac{1}{2}}:=\kappa_{\lambda,\eta}>0$.
	\qed
	
	\vskip2mm
	\noindent
	{\bf Lemma 2.5.} {\it If $|\lambda|<a\lambda_1$ and $\eta>\eta_0$, then $c=\inf_{u\in\mathcal{K}}I(u)>0$ can be achieved.
	}
	\vskip0mm
	\noindent
	{\bf Proof.}  Firstly, we claim that $c>0$. Indeed, for any $u\in\mathcal{K}$, according to Lemma 2.4, we have:
	\begin{eqnarray}\label{d1}
		I(u)=I(u)-\frac{1}{4}\langle I'(u),u\rangle=\frac{1}{4}I_{\lambda}(u)\ge \frac{1}{4}\left(a-\frac{|\lambda|}{\lambda_1}\right)\|u\|^2> \frac{1}{4}\left(a-\frac{|\lambda|}{\lambda_1}\right)\kappa_{\lambda,\eta}^2>0,
	\end{eqnarray}
	which implies that $c>0$.
	\par
	Next, we demonstrate that $c$ is achieved. Let $\{u_n\}\subset \mathcal{K}$ such that $\lim_{n\to\infty}I(u_n)=c$. Then from (\ref{d1}), it is straightforward to see that $\{u_n\}$ is bounded in $W^{1,2}_0(\Omega)$. Hence, since  $W^{1,2}_0(\Omega)$ is pre-compact, up to a subsequence,  there exists a $u_b\in W^{1,2}_0(\Omega)$ such that
	\begin{eqnarray}\label{d2}
		u_n\to u_b, \mbox{ as } n\to\infty,  \mbox{ in }W^{1,2}_0(\Omega).
	\end{eqnarray}
	Then Lemma 2.4 and (\ref{d2}) show that $\|u_b\|\ge \kappa_{\lambda,\eta}>0$. Hence, $u_b\not=0$.
\par
Moreover, by Lemma 2.1, for all $1\le \gamma< +\infty$, we have
	\begin{eqnarray}\label{d3}
		u_n\to u_b\;\;\text{in}\;L^{\gamma}(\Omega),\;u_n(x)\to u_b(x)\;\;\text{for all }\;x\in\Omega.
	\end{eqnarray}
Since $\{u_n\}\subset\mathcal{K}$, by $(\ref{d2})$ and $(\ref{d3})$, we get
	\begin{eqnarray}
		a\|u_b\|^2+b\|u_b\|^4
		& = & \lim_{n\to\infty}\left[a\|u_n\|^2+b\|u_n\|^4\right]\nonumber\\
		& = & \lim_{n\to\infty}\left[\lambda\int_{\Omega}|u_n|^2d\mu+\eta\int_{\Omega}|u_n|^4d\mu\right]\nonumber\\
		& = & \lambda\int_{\Omega}|u_b|^2d\mu+\eta\int_{\Omega}|u_b|^4d\mu,
	\end{eqnarray}
	which means $\langle I'(u_b),u_b\rangle= 0$. Hence, $u_b\in\mathcal{K}$. Furthermore, by (\ref{d2}) and the continuity of $I$, we have $I(u_b)=c$. The proof is completed. \qed

	\vskip2mm
	\noindent
	{\bf Lemma 2.6.} {\it If $\lambda<a\lambda_1$ and $\eta>\eta_0$, then for any given $u\in\mathcal{U}$, there exists a unique positive constant $s_u$ such that $s_uu\in\mathcal{K}$ and $I(s_u u)=\max_{s>0}I(su)$.
	}
	\vskip0mm
	\noindent
	{\bf Proof.}
	For any given $u\in\mathcal{U}$, it follows from  (\ref{3.2}) and (\ref{oo4}) that $su\in\mathcal{K}$ if and only if the positive constant $s$ satisfies
	\begin{eqnarray}\label{d4}
		\langle I'(su),su\rangle
		=s^2I_{\lambda}(u)+s^4I_{\eta}(u)
		=0.
	\end{eqnarray}
	Note that
	$
	I_{\eta}(u)<0.
	$
	Thus (\ref{d4}) and the fact $I_{\lambda}(u)>0$  show that
	\begin{eqnarray}\label{d5}
		s=s_u=\sqrt{-\frac{I_{\lambda}(u)}{I_{\eta}(u)}}>0.
	\end{eqnarray}
	which obviously is unique on $(0,+\infty)$.
	Furthermore, we define $J_u:(0,\infty)\to\R$ by $J_u(s)=I(su)$. Therefore, we can deduce that
	\begin{eqnarray}\label{d6}
		J'_u(s)	=\langle I'(su),u \rangle=\frac{1}{s}\langle I'(su),su \rangle=sI_{\lambda}(u)+s^3I_{\eta}(u),
	\end{eqnarray}
	which demonstrates that a positive constant $s$ is a critical point of $J_u$ if and only if $su\in\mathcal{K}$. Consequently, (\ref{d5}) shows that $s_u$ is a unique critical point of the function $J_u$.
	\par
	Next, we prove that  $J_u(s_u)=I(s_u u)=\max_{s>0}J_u(s)$.
	Since $u\in\mathcal{U}$, by (\ref{d4}) and (\ref{d6}), we get
	\begin{eqnarray}
		J''_u(s_u)=I_{\lambda}(u)+3s^2_uI_{\eta}(u)=2s^2_uI_{\eta}(u)<0.
	\end{eqnarray}
	Consequently,  $J_u(s_u)=\max_{s>0}J_u(s)$ and $s_u$ is the unique maximum point.
	\qed

	\vskip2mm
	\noindent
	{\bf Proof of Theorem 1.1.}\;  By Lemma 2.5, we have $I(u_b)=c$. We only need to prove that $I'(u_b)=0$.  Proving by contradiction, we suppose $I'(u_b)\not=0$.   Then, it follows from $I\in C^1(W_0^{1,2}(\Omega),\R)$ that  there are $\delta>0$ and $\theta>0$ such that $\|I'(v)\|_*\ge \theta$, for all $\|v-u_b\|\le 3\delta$, where $\|\cdot\|_*$ is the norm on the dual space of $W_0^{1,2}(\Omega)$.	
\par
 Since $u_b\in\mathcal{K}\subset \mathcal{U}$, $\langle I'(u_b),u_b\rangle=0$. Then Lemma 2.6 implies  that, for $s\in(0,1)\cup(1,+\infty)$,
	\begin{eqnarray}
		\label{c21}
		I(su_b)<I(u_b)=c.
	\end{eqnarray}
 Since $I_{\lambda}(u_b)>0$, by the continuity of $I_{\lambda}$, there exists $\sigma\in(0,1)$ small enough satisfying
 \begin{eqnarray}\label{mm1}
 \min_{s\in[1-\sigma,1+\sigma]}I_{\lambda}(su_b)>0.
\end{eqnarray}
 Denote $\Pi:=(1-\sigma,1+\sigma)$ and $\psi(s)=su_b$, $s\in\Pi$. Then by $(\ref{c21})$, we obtain that
	\begin{eqnarray}
		\label{c22}
		c_0:=\max_{\partial\Pi}I\circ\psi<c.
	\end{eqnarray}
	Let $\varepsilon:=\min\{(c-c_0)/3,\theta\delta/8\}$ and $S_{\delta}:=B(u_b,\delta)$, it follows from Lemma 2.3 in \cite{Willem2012} that there is $\gamma\in C([0,1]\times W^{1,2}_0(\Omega),W^{1,2}_0(\Omega))$ satisfying\\
	$(a)\;\;\gamma(\alpha,v)=v\;\;\text{if}\;\;\alpha=0\;\;\text{or}\;\;v\not\in I^{-1}([c-2\varepsilon,c+2\varepsilon]\cap S_{2\delta})$;\\
	$(b)\;\;I(\gamma(\alpha,v))<c\;\;\text{for all}\;\;v\in S_{\delta}\;\;\text{with}\;\;I(v)\le c\;\;\text{and}\;\;\alpha\in(0,1]$;\\
	$(c)\;\;I(\gamma(\alpha,v))\le I(v)\;\;\text{for all}\;\;v\in W^{1,2}_0(\Omega)\;\;\text{and}\;\;\alpha\in [0,1]$.
	\par
Next, we prove that
	\begin{eqnarray}
		\label{c23}
		\max_{s\in \Pi}I(\gamma(\alpha,\psi(s)))<c,\;\;\forall \alpha\in(0,1].
	\end{eqnarray}
	Indeed,  in view of $(b)$, it holds that $\max_{\{s\in\Pi:\psi(s)\in S_{\delta}\}}I(\gamma(\alpha,\psi(s)))<c,\;\;\forall \alpha\in (0,1]$. Moreover, it follows from $(c)$ and (\ref{c21}) that
	$$
	\max_{\{s\in\Pi:\psi(s)\not\in S_{\delta}\}}I(\gamma(\alpha,\psi(s)))\le \max_{\{s\in\Pi:\psi(s)\not\in S_{\delta}\}}I(\psi(s))<c,\;\;\forall \alpha\in [0,1].
	$$
	Thus, $(\ref{c23})$ can be proved.
	\par
	By $(a)$ and (\ref{mm1}), we have  $I\circ\gamma(0,su_b)=I(su_b)>0$ for all $s\in\Pi$.  Then it follows from the continuity  of $I\circ\gamma$ and the feature of guarantee code  that for any given  $s\in\Pi$, there exists $\alpha_s\in (0,1]$ such that
	\begin{eqnarray}
		\label{c24}
		I\circ\gamma(\alpha_s,su_b)>0.
	\end{eqnarray}
	Next, we prove that $\gamma((0,1],\psi(\Pi))\cap\mathcal{K}\not=\emptyset$. Let $\chi(s):=\gamma(\alpha_s,\psi(s))$,
	$
	\Psi_0(s):=\langle I'(\psi(s)),\psi(s)\rangle:=\varphi_{u_b}(s)
	$
and
	$\Psi_1(s):=\langle I'(\chi(s)),\chi(s)\rangle$.
	\par
	From $u_b\in\mathcal{K}$, (\ref{d4}) and (\ref{b3}), we can obtain that $\frac{\partial \varphi_{u_b}(s)}{\partial s}|_{s=1}=2I_{\lambda}(u_b)+4I_{\eta}(u_b)=-2I_{\lambda}(u_b)<0$. Since $\Psi_0(s)$ is continuously differentiable on $\Pi$ and $s=1$ is the unique isolated zero point of $\Psi_0$ on $\Pi$ by $u_b\in\mathcal{K}$, by using the Brouwer degree theory, we deduce that $\deg(\Psi_0,\Pi,0)=-1$. By $c_0<c-2\varepsilon$ and $(a)$, we get $\gamma(\alpha_s,\psi(s))=\psi(s),\;\forall s\in\partial\Pi,\;\alpha_s\in[0,1]$. So,  $\Psi_0(s)=\Psi_1(s)$ on $\partial\Pi$. By the homotopy invariance of Brouwer degree, we have $\deg(\Psi_1,\Pi,0)=-1$, which shows that $\Psi_1(s_0)=0$ for some $s_0\in\Pi$. According to $(\ref{c24})$, we obtain $\chi(s_0)=\gamma(\alpha_{s_0},\psi(s_0))\not= 0$. So $ \chi(s_0)\in\mathcal{K}$. Thus, $\gamma((0,1],\psi(\Pi))\cap\mathcal{K}\not=\emptyset$. Thanks to $(\ref{c23})$, we conclude a contradiction. Thus, we finish the proof.
	\qed
	\vskip2mm
	\par
	Finally, in this section, we present a minimax characterization of the critical value $c$.
	\vskip2mm
	\noindent
	{\bf Lemma 2.7.} {\it There exists the following minimax characterization:}
	\begin{eqnarray}
		\label{n3}
		\inf_{u\in\mathcal{K}}I(u)
		= c
		= \inf_{u\in \mathcal{U}}\max_{s>0}I(su).
	\end{eqnarray}
	\vskip0mm
	\noindent
	{\bf Proof.} For $u\in\mathcal{K}$, by Lemma 2.6, we can obtain that $s_{u}= 1$ and
	\begin{eqnarray} \label{n2}
		I(u)=\max_{s> 0}I(su).
	\end{eqnarray}
	Then it follows from (\ref{n2}) and $\mathcal{K} \subset\mathcal{N}$ that
	$$
	\inf_{u\in\mathcal{K}}I(u)
	=  \inf_{u\in \mathcal{K}}\max_{s>0}I(su)
	\geq \inf_{u\in \mathcal{U}}\max_{s> 0}I(su).
	$$
	Moreover, by Lemma 2.6, for any $u\in \mathcal{U}$, there exists $s_u>0$ such that $s_u u\in \mathcal{K}\subset\mathcal{U}$. Then there exists
	\begin{eqnarray*}
		\label{n1}
		\max_{s> 0} I(su)
		\geq I(s_u u)
		\geq \inf_{u\in\mathcal{U}}I(u).
	\end{eqnarray*}
	Furthermore, we have
	$$
	\inf_{u\in \mathcal{U}}\max_{s> 0} I(su)
	\geq \inf_{u\in\mathcal{N}}I(u).
	$$
	Hence, (\ref{n3}) holds.

	\vskip2mm
	{\section{The least energy solution for $(\ref{t2})$}}
	\setcounter{equation}{0}
	
	In this section, we prove the existence of the least energy  solution for $(\ref{t2})$ and complete the proof of Theorem 1.2. There exist some differences between
the proofs of Theorem 1.1 and Theorem 1.2, whose reason is that  $W_h^{1,2}(V)$ is of infinite dimension but  $W^{1,2}_0(\Omega)$ is of finite dimension so that the boundness of $\{u_n\}\subset \mathcal{M}$ just leads to the existence of a  weakly convergent subsequence rather than a strong convergent subsequence.
 We shall work on $W_h^{1,2}(V)$ and need the following Sobolev embedding theorem.
  \vskip2mm
	\noindent
	{\bf Lemma 3.1.} (\cite{Yang}, $p=2$) {\it If $\mu(x)\geq \mu_0>0$ and $(H_1)$ holds. Then  $W_h^{1,2}(V)$ is continuously embedded into $L^r(V)$ for all $2\le r\le \infty$, and the following inequalities hold:
		\begin{eqnarray}\label{h5}
			\|u\|_{\infty}\le \frac{1}{h_0^{\frac{1}{2}}\mu_0^{\frac{1}{2}}} \|u\|_{W_h^{1,2}(V)}
		\end{eqnarray}
		and
		\begin{eqnarray} \label{h6}
			\|u\|_{V,r}\le \mu_0^{\frac{2-r}{2r}}h_0^{-\frac{1}{2}}\|u\|_{W_h^{1,2}(V)}\ \ \mbox{for all } 2\le r< \infty.
		\end{eqnarray}
	Furthermore, if $(H_2)$ also holds, then  $W_h^{1,2}(V)$ is compactly embedded into $L^r(V)$ for all $2\le r\le \infty$. }
	
	\vskip2mm
	\par
	Let $J_{\lambda}(u):=a\|u\|_{ W^{1,2}_h(V)}^2-\lambda\int_{V}|u|^2d\mu$ and $J_{\eta}(u):=b\|u\|_{ W^{1,2}_h(V)}^4-\eta\int_{V}|u|^4d\mu$. Then the energy functional given by (\ref{cm2}) can be written as
	\begin{eqnarray}\label{m2}
		J(u)=\frac{1}{2}J_{\lambda}(u)+\frac{1}{4}J_{\eta}(u),\;\;\;u\in W^{1,2}_h(V).
	\end{eqnarray}
Moreover, a standard argument implies that $J(u)$ belongs to $C^1(W^{1,2}_h(V),\R)$, and
	\begin{eqnarray}\label{4.2}
		\langle J'(u),v\rangle=a\int_{V}(\nabla u\nabla v+h(x)uv)d\mu-\lambda\int_{V}uvd\mu+b\|u\|^2_{W^{1,2}_h(V)}\int_{V}(\nabla u\nabla v+h(x)uv)d\mu-\eta\int_{V}|u|^2uvd\mu
	\end{eqnarray}
	for any $u,v\in W^{1,2}_h(V)$ (for example, see \cite{Yang}). Thus, the critical points of $J$ correspond to the weak solutions of (\ref{t2}), which actually are the point-wise solutions of (\ref{t2}) (for example, see Proposition 3.2 in \cite{Yang}).
	
	\par
	We denote the following constraint by $\mathcal{M}$:
	\begin{eqnarray}
		\label{m1}
		\mathcal{M}=\{u\in W^{1,2}_h(V),u\not=0:\langle J'(u),u\rangle=0\},
	\end{eqnarray}
	and define the following minimization problem
	$$
	d=\inf_{u\in \mathcal{M}}J(u).
	$$
	If $u\in W^{1,2}_h(V)$ is a nontrivial solution of $(\ref{t2})$, then by (\ref{4.2}), it holds that
	\begin{eqnarray}
		\label{m3}
		0=\langle J'(u),u\rangle=J_{\lambda}(u)+J_{\eta}(u).
	\end{eqnarray}
	Thus, we deduce from $\lambda<a\lambda_1^*$ and $(\ref{m3})$ that
	\begin{eqnarray*}
		J_{\eta}(u)=b\|u\|^4_{W^{1,2}_h(V)}-\eta\int_{V}|u|^4d\mu= -J_{\lambda}(u)
		=-a\|u\|^2_{W^{1,2}_h(V)}+\lambda\int_{V}|u|^2d\mu<-a\|u\|^2_{W^{1,2}_h(V)}+a\lambda_1^*\int_{V}|u|^2d\mu\le0,
	\end{eqnarray*}
	that is
	\begin{eqnarray*}
		\eta>\frac{b\|u\|^4_{W^{1,2}_h(V)}}{\int_{V}|u|^4d\mu},\;\;\;\text{for any }u\in\mathcal{M}.
	\end{eqnarray*}

 \par
 The proofs of Lemma 3.2-3.5 below are almost the same as Lemma 2.2-2.5. We omit the details.
	\vskip2mm
	\noindent
	{\bf Lemma 3.2.} {\it Suppose $\eta\le\eta_0^*=\inf_{u\in W^{1,2}_h(V)\setminus\{0\}}\frac{b\|u\|^4_{W^{1,2}_h(V)}}{\int_{V}|u|^4d\mu}$. Then $(\ref{t2})$ has no any nontrivial solution.}
	\vskip2mm
	\par
	We denote
	\begin{eqnarray*}
		\label{o14}
		\mathcal{N}=\Bigg\{u\in W^{1,2}_h(V),u\not=0:J_{\eta}(u)<0\Bigg\}
		=\left\{u\in W^{1,2}_h(V),u\not=0:b\|u\|^4_{W^{1,2}_h(V)}-\eta\int_{V}|u|^4d\mu<0\right\}.
	\end{eqnarray*}

	\vskip2mm
	\noindent
	{\bf Lemma 3.3.} {\it If $\lambda<a\lambda_1^*$ and $\eta>\eta_0^*$, then $\mathcal{N}\not=\emptyset$ and $\mathcal{M}\subset\mathcal{N}$.}
	
	\vskip2mm
	\noindent
	{\bf Lemma 3.4.} {\it If $|\lambda|<a\lambda_1^*$ and $\eta>\eta_0^*$, then $\|u\|_{W^{1,2}_h(V)}>\kappa_{\lambda,\eta}^*:=\left(\frac{a-\frac{\lambda}{\lambda_1^*}}{\eta^*}\mu_0h_0^2\right)^{\frac{1}{2}}$ for every $u\in\mathcal{M}$.
	}

	\vskip2mm
	\noindent
	{\bf Lemma 3.5.} {\it If $\lambda<a\lambda_1^*$ and $\eta>\eta_0^*$, then for any given $u\in\mathcal{N}$, there exists a unique positive constant $s_u$ such that $s_uu\in\mathcal{M}$ and $J(s_u u)=\max_{s>0}J(su)$.
	}
\vskip2mm
	\par
     Since $W_h^{1,2}(V)$ is of infinite dimension,  the boundness of $\{u_n\}\subset \mathcal{M}$ can not lead to the existence of  a strongly convergent subsequence as done in Theorem 1.1 to prove $d=\inf_{u\in\mathcal{M}}J(u)$ to be achieved. In order to overcome such difficulty, we need the following lemmas.
\vskip2mm
	\noindent
	{\bf Lemma 3.6.} {\it For any given $u\in\mathcal{N}$, if $\langle J'(u),u\rangle\le 0$, then $0<s_u\le 1$.
	}
	\vskip0mm
	\noindent
	{\bf Proof.}
	For any given $u\in\mathcal{N}$, it follows from  $(\ref{m1})$ that $su\in\mathcal{M}(s>0)$ if and only if  $s$ satisfies
	\begin{eqnarray}\label{o4}
		\langle J'(su),su\rangle
		=s^2J_{\lambda}(u)+s^4J_{\eta}(u)
		=0.
	\end{eqnarray}
 Since $s_uu\in\mathcal{M}$, by (\ref{o4}), we derive
	\begin{eqnarray}\label{o7}
		s^2_uJ_{\eta}(u)=-J_{\lambda}(u).
	\end{eqnarray}
	Furthermore, by  $\langle  J'(u),u\rangle\le 0$ and $J_{\lambda}(u)>0$, we have
	\begin{eqnarray}
		\label{o8}
		J_{\eta}(u)\le -J_{\lambda}(u)<0.
	\end{eqnarray}
	By (\ref{o7}) and (\ref{o8}),  we deduce that
	\begin{eqnarray*}
		-s^2_uJ_{\lambda}(u)\ge s^2_uJ_{\eta}(u)=-J_{\lambda}(u).
	\end{eqnarray*}
	Since $J_{\lambda}(u)>0$,  we conclude $0<s_u\le 1$.
	\qed

	\vskip2mm
	\noindent
	{\bf Lemma 3.7.} {\it If $|\lambda|<a\lambda_1^*$ and $\eta>\eta_0^*$, then $d=\inf_{u\in\mathcal{M}}J(u)>0$ is achieved.
	}
	\vskip0mm
	\noindent
	{\bf Proof.}  Firstly, we claim that $d>0$. Indeed, for any $u\in\mathcal{M}$, according to Lemma 3.4, we have
	\begin{eqnarray}\label{o1}
		J(u)=J(u)-\frac{1}{4}\langle J'(u),u\rangle=\frac{1}{4}J_{\lambda}(u)\ge \frac{1}{4}\left(a-\frac{|\lambda|}{\lambda_1^*}\right)\|u\|^2_{W^{1,2}_h(V)}> \frac{1}{4}\left(a-\frac{|\lambda|}{\lambda_1^*}\right)(\kappa_{\lambda,\eta}^*)^2>0,
	\end{eqnarray}
	which implies that $d>0$.
	\par
	Next, we demonstrate that $d$ is achieved. Let $\{u_n\}\subset \mathcal{M}$ that $\lim_{n\to\infty}J(u_n)=d$. Then from (\ref{o1}), it is straightforward to see that $\{u_n\}$ is bounded in $W^{1,2}_h(V)$. Hence, there exists $u_b^*\in W^{1,2}_h(V)$ such that
	\begin{eqnarray}\label{o2}
		u_n\rightharpoonup u_b^*, \mbox{ as } n\to\infty.
	\end{eqnarray}
	Moreover, by Lemma 3.1, for all $2\le p\le +\infty$, we have
	\begin{eqnarray}\label{o3}
		u_n\to u_b^*\;\;\text{in}\;L^p(V),\;u_n(x)\to u_b^*(x)\;\;\text{for all }\;x\in V.
	\end{eqnarray}
	Lemma 3.4, the weak lower semi-continuity of norm and (\ref{o2}) show that $\|u_b^*\|_{W^{1,2}_h(V)}\ge \kappa_{\lambda,\eta}^*>0$. Hence, $u_b^*\not=0$.
	\par
	On the other hand, since $\{u_n\}\subset\mathcal{M}$, by  the weak lower semi-continuity of norm, $(\ref{o2})$ and $(\ref{o3})$, we get
	\begin{eqnarray*}
		a\|u_b^*\|^2_{W^{1,2}_h(V)}+b\|u_b^*\|^4_{W^{1,2}_h(V)}
		& \le & \lim_{n\to\infty}\left[a\|u_n\|^2_{W^{1,2}_h(V)}+b\|u_n\|^4_{W^{1,2}_h(V)}\right]\nonumber\\
		& = & \lim_{n\to\infty}\left[\lambda\int_{V}|u_n|^2d\mu+\eta\int_{V}|u_n|^4d\mu\right]\nonumber\\
		& = & \lambda\int_{V}|u_b^*|^2d\mu+\eta\int_{V}|u_b^*|^4d\mu,
	\end{eqnarray*}
	which means $\langle J'(u_b^*),u_b^*\rangle\le 0$.
According to Lemma 3.5 and Lemma 3.6, there exists $0<s_{u_{b^*}}\le 1$ such that $u'_b:=s_{u_b}u_b^*\in\mathcal{M}$. Thanks to $(\ref{o3})$, we have
\begin{eqnarray*}
	d&\le& J(u'_b)=J(s_{u_b}u_b^*)=J(s_{u_b}u_b^*)-\frac{1}{4}\langle J(s_{u_b}u_b^*),s_{u_b}u_b^*\rangle\\
	&=&\frac{s^2_{u_{b^*}}}{4}J_{\lambda}(u_b^*)\le\frac{1}{4}J_{\lambda}(u_b^*)=\frac{a}{4}\|u_b^*\|^2_{W^{1,2}_h(V)}-\frac{\lambda}{4}\int_{V}|u_b^*|^2d\mu\\
	&\le& \liminf_{n\to\infty}\left[\frac{a}{4}\|u_n\|^2_{W^{1,2}_h(V)}-\frac{\lambda}{4}\int_{V}|u_n|^2d\mu\right]\\
	&=&\liminf_{n\to\infty}\left[J(u_n)-\frac{1}{4}\langle J'(u_n),u_n\rangle\right]=d.
\end{eqnarray*}
Therefore, $s_{u_{b^*}}=1$, and $d$ is achieved by $u'_b=u_b^*\in\mathcal{M}$. The proof is completed.
\qed

	\vskip2mm
	\noindent
	{\bf Proof of Theorem 1.2.}\; By Lemma 3.7, we have $J(u_b^*)=d$. We only need to prove that $J'(u_b^*)=0$. In fact, with the help of Lemma 3.5, the proof is the same as Theorem 1.1 with replacing the functional $I$ and $u_b$ with $J$ and $u_b^*$, respectively. We omit the details. \qed
	
	\vskip2mm
	\par
	Finally, using the same proof as Lemma 2.7, we also present a minimax characterization of the critical value $d$.
	\vskip2mm
	\noindent
	{\bf Lemma 3.8.} {\it There exists the following minimax characterization:}
	\begin{eqnarray*}
		\label{oo3}
		\inf_{u\in\mathcal{M}}J(u)
		= d
		= \inf_{u\in \mathcal{N}}\max_{s>0}J(su).
	\end{eqnarray*}

	\vskip2mm
	\noindent
	{\bf Acknowledgments}\\
	This work is supported by Yunnan Fundamental Research Projects of China (grant No: 202301AT070465).

\vskip2mm
\renewcommand\refname{References}
{}

\end{document}